\documentclass[]{article}
\usepackage{fancyhdr, epsfig, times}
\usepackage[amsthm]{ntheorem}
\usepackage{amsmath}
\usepackage{amssymb}
\usepackage{multicol}
\usepackage{graphicx}

\theoremstyle{plain}\newtheorem{owndefinition}{Definition}
\theoremstyle{plain}
\theoremstyle{definition}
\theoremstyle{remark}\newtheorem*{ownremark}{Remark}

\begin{document}

\newpage
\thispagestyle{empty}

\newpage

\title{Approximating the Euclidean Circle in the Square Grid using Neighbourhood Sequences}

\author{J\'anos Farkas$^1$, Szabolcs Baj\'ak$^2$, Benedek Nagy$^{1,3}$
\thanks{The research is partly supported by grants F043090 and T049409 of the National 
 Foundation of Scientific Research of Hungary.}\\
{\small $^1$ Faculty of Informatics,}\\
{\small University of Debrecen, H-4032 Debrecen, Egyetem t\'er 1, Hungary.}\\
{\small $^2$ Institute of Mathematics,}\\
{\small University of Debrecen, H-4032 Debrecen, Egyetem t\'er 1, Hungary.}\\
{\small $^3$ Research Group on Mathematical Linguistics,}\\
{\small Rovira i Virgili University, Tarragona, Spain.}\\
{\footnotesize \tt \{farkasj|nbenedek\}\mbox{@}inf.unideb.hu,
bajaksz\mbox{@}delfin.klte.hu}}
\date{}
\maketitle

\pagenumbering{arabic}

\begin{abstract}
\par Distance measuring is a very important task in digital geometry and di\-gi\-tal image processing. Due to our natural approach to geometry we think of the set of points that are equally far from a given point as a Euclidean circle. Using the classical neighbourhood relations on digital grids, we get circles that greatly differ from the Euclidean circle.
\par In this paper we examine different methods of approximating the Euclidean circle in the square grid, considering the possible motivations as well. We compare the pe\-ri\-me\-ter-, area-, curve- and noncompactness-based  approximations and examine their realization using neighbourhood sequen\-ces.
\par We also provide a table which summarizes our results, and can be used when developing applications that support neighbourhood sequences.
\end{abstract}

\section{Introduction}
\par The classical digital geometry started with \cite{rosenf}, where the authors defined the two basic 
neighbourhood relations on the square grid. The topic is well developed due to 
people of image processing and computer graphics communities. In \cite{1} the authors used the so-called neighbourhood
sequences to vary 
the neighbourhood criterion in a path. They used only periodic neighbourhood sequences in their a\-na\-ly\-sis.
Some properties of distances based on neighbourhood sequences
are detailed in \cite{2o}. The concept of neighbourhood sequences was extended to not necessarily periodic
sequences in \cite{3}. These general neighbourhood sequences were analysed in \cite{6,pub},
formulae to compute distances are presented in \cite{isp}.
 One of the main problems of digital geometry is the approximation of the Euclidean circle (for instance in \cite{2}). It is the topic of the present paper as well. \\
After the formal definitions (Section 2) some previous
results on the topic are recalled (Section 3). The approximations of the circle by neighbourhood sequences (i.e. the digital discs) are octagons. The descriptors of these octagons are presented in Section 4.
In Section 5 several approaches of approximation are detailed,
in Section 6 a detailed example is shown, while in Section 7 we summarize the
results. Finally a conclusion closes the paper.

\section{Definitions} 
\par Our aim is to provide the best approximation of a circle with a given radius using neighbourhood sequences. In order to formulate our results, we have to recall the following definitions \cite{1,3,2}. 
Because a circle is a planar shape, we restrict our considerations to two dimensions.

\begin{owndefinition} 
    Let $p$ and $q$ be two points in $\mathbb{Z}^2$ and $j \in \{1, 2\}$. The $i$th coordinate of the point $p$ is indicated by
    $Pr_i(p)$. The points $p$ and $q$ are \emph{$j$-neighbours in two dimensions}
    if the following two conditions hold:
    \begin{itemize}
        \item $|Pr_i(p)-Pr_i(q)| \leq 1 $ \quad $(\forall i \in \{1, 2\})$,
        \item $\sum_{i=1}^2 |Pr_i(p)-Pr_i(q)| \leq j$.
    \end{itemize}
\end{owndefinition}

\par $1$-neighbourhood corresponds to the classical $4$-neighbourhood and $2$-neigh\-bour\-hood to $8$-neighbourhood, as it can be seen on Fig. 1 left and right, respectively. In two dimensions, $j$-neighbourhood means that we can step along at most $j$ Cartesian axes to reach a neighbour.\\

\begin{figure}[ht]
  \centerline{\epsfig{file=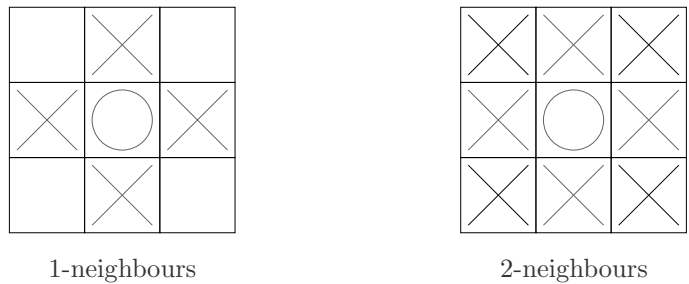} }
    \caption{The X-marked points are $j$-neighbours of the O-marked one} \label{fig1}
\end{figure}


\par In a neighbourhood sequence one can give a sequence of neighbourhoods, in which a $j$-neighbourhood
is represented by the number $j$. The $i$th element of such a sequence defines the neighbourhood we can
use as we take the $i$th step going further and further from the starting point.
Formally in the $n$ dimensional digital space the neighbourhood sequences are defined in the following way:

\begin{owndefinition} 
    The infinite sequence $B=(b_i)_{i=1}^{\infty}$, where $b_i \in \{1, \ldots, n\}$ for all $i \in \mathbb{N}$, is called
    a \emph{generalized $n$D-neighbourhood sequence}.
\end{owndefinition}

In $\mathbb{Z}^2$, neighbourhood sequences are infinite sequences of values $1$ and $2$. Now we can define the path leading from one point to another, its length, and the distance between two points.

\begin{owndefinition} 
    Let $p$ and $q$ be two points in $\mathbb{Z}^2$ and $B$ a generalized $2$D-neighbourhood sequence. The point
    sequence $\Pi(p,q;B)$ - which has the form $p=p_0,p_1,\dots,p_m = q$, where $p_{i-1}$ and $p_i$ are
    $b_i$-neighbours for $i \in \{1, 2, 3, \ldots\}$ - is called a \emph{path from $p$ to $q$ determined by $B$}. The
    \emph{length $|\Pi(p,q;B)|$ of the path} $\Pi(p,q;B)$ is $m$.
\end{owndefinition}

\begin{owndefinition} 
    Let $p$ and $q$ be two points in $\mathbb{Z}^2$ and $B$ a generalized $2$D-neighbourhood sequence.
    The/a
    shortest path from $p$ to $q$ is denoted by $\Pi^*(p,q;B)$. The \emph{distance between $p$ and $q$}
    is defined as the length of the minimal path, and is written as
    \[
        d(p,q;B)=|\Pi^*(p,q;B)|.
    \]
\end{owndefinition}

\par
\begin{figure}[ht]
  \centerline{\epsfig{file=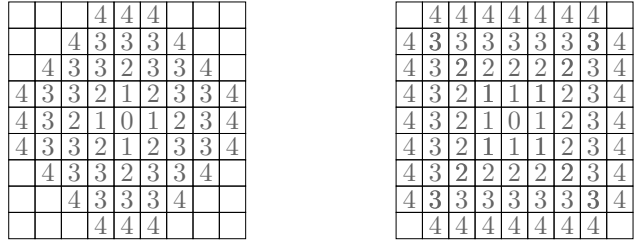} }
%
\caption{Examples of digital discs with a radius of 4}\label{f2}
\end{figure}

 As one can see, neighbourhood sequences occupy digital octagons after every step (see Fig. 2).
This means, that if we use neighbourhood sequences on the square
grid, the only shape we can use to approximate circles are special
octagons, having $135^\circ$ inner angles \cite{2}, or squares with right angles in degenerated cases.
The octagon is degenerated if only 1 type of neighbourhood is used to generate it. On Fig. \ref{f2}, we
can see two discs with a radius of 4, generated by the
neighbourhood sequences 
$B=(1,1,2,1,...)$ and $B=(2,2,2,1)$, respectively. On the figures, each point of the digital discs is
indicated by a number, which is the distance of the point from the center point with the given sequence.


\smallskip

\section{Previous results}
\par In \cite{2} Hajdu and Nagy  used the isoperimetric ratio (\cite{1}) to approximate the Euclidean circle.

\begin{owndefinition} 
    We call
    \[ \kappa := \frac{P^2}{A}, \]
    the \emph{isoperimetric ratio} of the shape, where $P$ is the perimeter, and $A$ is the area of the shape.
\end{owndefinition}

\begin{ownremark}
    In \cite{1} and \cite{2} the isoperimetric ratio is referred as \emph{noncompactness ratio}.
\end{ownremark}

\par In Euclidean geometry, according to the isoperimetric inequality, the isoperimetric ratio is never less than $4\pi$, and is minimal for the circle. Let $k_1$ denote the number of value $1$ in the neighbourhood sequence $B$, and $k_2$ the number of value $2$. The number of steps we take in the neighbourhood sequence $B$ will be denoted by $k$, and $k=k_1+k_2$. Let $c$ denote the proportional frequency of value 2 in the sequence, thus $ c:=\frac{k_2}{k}=\frac{k_2}{k_1+k_2}$. The isoperimetric ratio for the convex hull of the possible octagons can be written in the form
\[ \kappa(c)=16\left( 1-2\left(2-\sqrt{2}\right)\frac{c(1-c)}{c(2-c)+1}\right). \]
\par It turns out that $\kappa(c)$ is minimal at $c = \sqrt{2}-1$. If $c$ satisfies this equality, the convex hull of the pixel centers becomes a regular octagon, thus the problem of approximating the circle turns out to be the problem of approximating the regular octagon. Note that because of the irrationality of $\sqrt{2}-1$, $c$ can only be equal to this in the ideal case of $k \rightarrow \infty$. For any $\varepsilon > 0$ there is
a value of $k_0$ such that for any $k > k_0$ \ $\exists k_2 : \frac{k_2}{k} - \left( \sqrt{2} -1  \right) < \varepsilon$.
\par Let us remark that the use of the isoperimetric ratio can be avoided. Taking into consideration the connection between the isoperimetric ratio and the isoperimetric problem, we can say that Hajdu and Nagy solved the isoperimetric problem for the corresponding octagons in their approximations on the square grid. The proof is much easier if we solve the isoperimetric problem directly for these octagons. Further in this article we take into account the circle's isoperimetric property by restricting our examinations only to regular octagons.
This can only be achieved if we assume that the approximating octagon is large enough to have more than one value 1 and 2 in its sequence $(k_1, k_2 \geq 1)$.
\par Hajdu and Nagy also gave the sequence which in every step generates the octagon closest to the regular octagon. In this paper we use a different approach. We assume that we have a circle with a particular radius $r$. We want to determine the sequence that leads to the best approximating regular octagon. We can also approximate circles with positive real radius, which could not be done by using the old approach. We use different measures to describe the generated octagons, which we will call \emph{`descriptors'}. We also investigate several possible definitions for `approximation'.

\section{Descriptors of the octagons}
\par The sequence-generated octagons can be easily described by their sidelengths $a$ and $b$ (see Fig. 3). We use three types of measure to give the length of the sides: the pixel, the inner convex hull (or inner octagon) and the outer convex hull (or outer octagon) based descriptors.

\medskip

\par \emph{Pixel based descriptors:} $a$ and $b$ are given in pixels. As it can be seen on Fig. 3, we assume that corner pixels belong to the horizontal sides ($a$).
\par \emph{Inner octagon:} in this case, $a$ and $b$ are equal to the corresponding sides of the convex hull of the centers of the pi\-xels.
\par \emph{Outer octagon:} the same as the inner octagon, except that we use the convex hull of the pixels as squares.

\begin{figure}[ht]
  \centerline{\epsfig{file=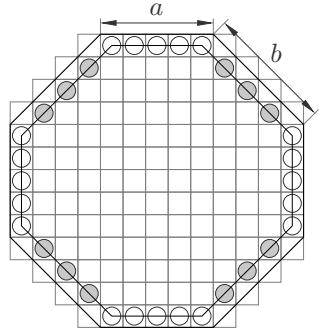} }
\caption{Descriptors}\label{f3}
\end{figure}


\par Table 1 shows the values for the sides, the perimeter and the area of the octagons.\\
\begin{center}
\begin{tabular}{c  c  c  c}
    \vspace{0.2cm}
    &\textbf{pixel} $(k_1\geq 1)$ & \textbf{inner hull} & \textbf{outer hull} \\ \vspace{0.2cm}
    $\mathbf{a}$ & $2k_2+1$ & $2k_2$        & $2k_2+1$\\ \vspace{0.2cm}
    $\mathbf{b}$ & $k_1-1$  & $\sqrt{2}k_1$ & $\sqrt{2}k_1$\\ \vspace{0.2cm}
    $\mathbf{P}$ & $4(a+b)$ & $4(a+b)$      & $4(a+b)$\\ \vspace{0.2cm}
    $\mathbf{A}$ & $(2k_2+1)^2 + 2k_1(4k_2+k_1+1)$ & $a^2+2\sqrt{2}ab+b^2$ & $a^2+2\sqrt{2}ab+b^2$
\end{tabular}\\
\end{center}
\begin{center} \emph{Table 1. Attributes of the octagons} \end{center}

\medskip
\par In the following section the inner and outer octagons give almost the same results (and they are equal to each other in the ideal case $k\rightarrow \infty$), so we omit the consideration of the outer octagon.
\par If our aim is only to make the octagon of the sequence the best approximating one in the sense that it is isoperimetric (it is as regular as it can be after a large fixed $k$ number of steps), we only need to solve the equation $a = b$ for each descriptors in the case of $k\rightarrow \infty$. Let $p$ denote the probability of the occurrence of value $2$ in the neighbourhood sequence, i.e.,
\[ p:=\lim_{k\to \infty} \frac{k_2}{k}. \]
\par Using this notation, we get that the regular octagons for large distances ($k\rightarrow \infty$) occur at $p = \frac{1}{3}$ when using the pixel-based descriptors and at $p = \sqrt{2}-1$ for the inner and outer convex hulls, just as in \cite{2}.
\par In \cite{2}, the concept of `convex hull' is equivalent to the descriptor we call inner octagon.

\section{Possible ways of approximating the circle}
\par In this section we split our ideas by answering the question: in what sense would we like to approximate the circle? We can use different measures to order the approximating octagons from `badly' approximating octagons (far from the result) to `well' approximating ones (close to the result). In this terminology we would like to minimize a distance function between the circle and the approximating shapes.
\par By inverting the results in Table 1 and using $a = b$, we can give the number of values 1 and 2 as a function of the sidelength $a$ (see Table 2). This way we also get the minimal length $k = k_1 + k_2$ of the sequence belonging to the best approximating regular octagon. In Table 2 $[.]$ denotes the rounding function.\\
\begin{center}
\begin{tabular}{c  c  c  c}
    \vspace{0.2cm}
    & \textbf{pixel $(k_1 \geq 1)$} & \textbf{inner octagon} & \textbf{outer octagon} \\ \vspace{0.2cm}
    $\mathbf{k_1}$ &$\left[a+1\right]$ & $\left[ \frac{a}{\sqrt{2}} \right]$ & $\left[ \frac{a}{\sqrt{2}} \right]$\\ \vspace{0.2cm}
    $\mathbf{k_2}$ & $\left[ \frac{a-1}{2} \right]$ & $\left[ \frac{a}{2} \right]$ & $\left[ \frac{a-1}{2} \right]$\\ \vspace{0.2cm}
    $\mathbf{p}$ & $\frac{1}{3}$ & $\sqrt{2}-1$ & $\sqrt{2}-1$\\
\end{tabular}
\end{center}
\begin{center} \emph{Table 2. Number of values 1 and 2 as a function of sidelength $a$} \end{center}
\par This means that we only need to compute $a$ in order to gain the appropriate sequence. In the following subsections we only provide the formulas from which we formulated our results.

\subsection*{Perimeter based approximation}
\par We would like to get the octagon having the same perimeter as the given circle. Formally we have to solve the following equation for $a_{perim}$:
\[  P_{regoct} = 8a = 2\pi r \]
\[ a = \frac{\pi}{4}\,r \]
\par We get different results using different descriptors. For each approach we need to get the ideal (continuous) length of $a_{ideal}$ (in this case $a_{perim}$), then solve the equations $a_{ideal} = a$ and $a_{ideal} = b$ by substituting the formulas of $a$ and $b$ from Table 1. This way we get the sequence of the best approximating regular octagon by having the values of $k_1$ and $k_2$. We follow the same method in all cases.

\subsection*{Area based approximation}
Let us determine the regular octagon with the same area as the given circle by solving the following equation for $a$:
\[ A_{regoct}=2\left(1+\sqrt{2}\right)a^2=\pi r^2 \]
\[ a = \sqrt{\frac{\pi}{2\left(1+\sqrt{2}\right)}}\,r \]

\subsection*{Inscribed circle based approximation}
\par The radius based approximations use the concept of the radius of the regular octagon, which is the distance between the center of the octagon and the sides, which is the same as the radius of the inscribed circle of the octagon. There are two types of radius based approximation: the inscribed circle based and the covering circle based method.
\par In the inscribed circle based approximation we would like to get the sequence that generates the regular octagon which is closest to the octagon having the given circle as its inscribed circle.

\par Since $r$ denotes the given radius, in this case it is the same as the radius of the generated octagon. We have to solve the following equation (for notions see Fig. 4):
\[ \tan \frac{\pi}{8} = \sqrt{2}-1 = \frac{a}{2r} \]
\[ a = 2\left(\sqrt{2}-1\right)\,r \]

\begin{figure}[ht]
  \centerline{\epsfig{file=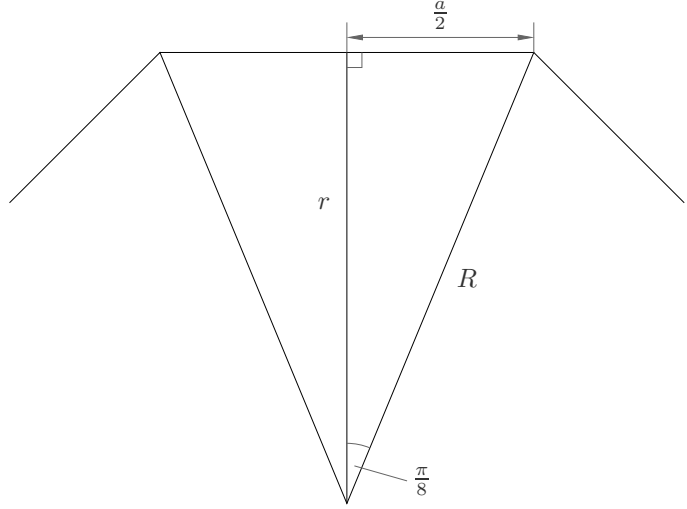} }
  \caption{Radius based approximation}\label{f4}
\end{figure}

\subsection*{Covering circle based approximation}
\par We need the sequence generating the octagon having the given circle as its co\-ve\-ring circle. This means that the radius of the given circle on Fig. 4 is denoted by $R$.
\[ R^2=\left( \frac{a}{2} \right)^2 + r^2 \]
\[ a = \sqrt{2-\sqrt{2}}\,R \]

\subsection*{Least squares difference}
\par Finally, we discuss curve based approximations. In these methods we use curve fitting, i.e., we search the regular octagon having its curve the closest to the curve of the given circle. By defining the distance between the two curves differently, we get two methods: the least square and the least sum of distances approximations. Due to symmetric reasons we only need to examine a fraction (namely one sixteenth) of the curves. To make the integration easier we convert the problem into a planar polar coordinate system (see Fig. 5), in which we denote the curve of the circle by $z_0(\alpha)$ and the curve of the octagon by $z(\alpha)$:
\[ z_0(\alpha)=r \qquad \qquad z(\alpha)=\frac{m}{\cos \alpha}. \]
\par  Using the least square method the distance between the two curves at $\alpha$ becomes
\[ f(\alpha,r):=\left(z(\alpha)-z_0(\alpha)\right)^2 .\]
The distance of the two curves can be determined by integrating these elementary distances:
\[ F(r):=\int\limits_{0}^{\frac{\pi}{8}} f(\alpha,r) d\alpha .\]

\par This way the search for the closest octagon becomes an extremum problem for $F(r)$. By solving $F'(r) = 0$ we get the optimal sidelength $a$, which is
\[ a = \frac{\pi}{4\left(\sqrt{2}+1\right)\ln \tan \left( \frac{5\pi}{16}\right)}\,r. \]


\par
\begin{figure}[ht]
  {\epsfig{file=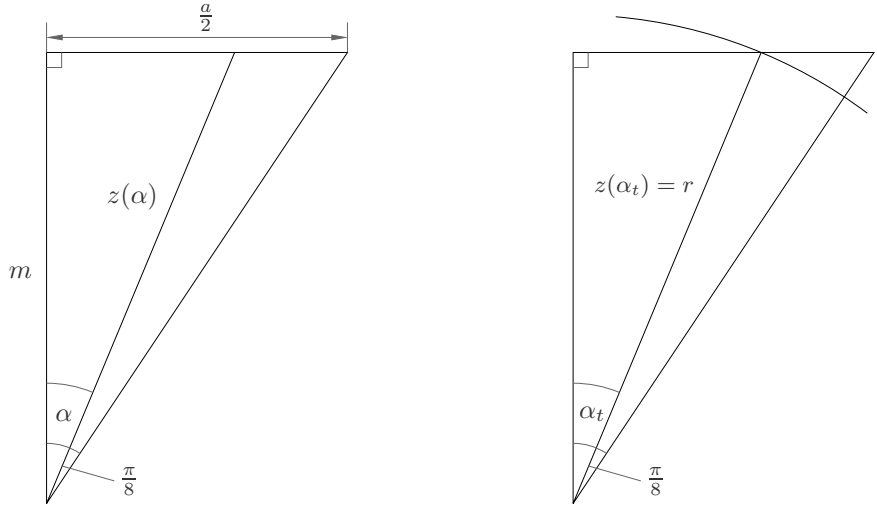} }
\caption{Curve based approximation}
\end{figure}

\subsection*{Least sum of distances}
\par In this case, the elementary distance is not squared, thus we have to manually take into account, that between the threshold angle $\alpha_t$ and $\frac{\pi}{8}$ the distance is negative.
\[ f(\alpha,r):=z_0(\alpha)-z(\alpha) \]
\[ z_0(\alpha_t)=z(\alpha_t) \quad \Rightarrow \quad
     \alpha_t=\arccos\left(\frac{\sqrt{2}+1}{2}\frac{a}{r}\right)   \]
\[ F(r):=\int\limits_{0}^{\alpha_t} f(\alpha,r) d\alpha - \int\limits_{\alpha_t}^{\frac{\pi}{8}} f(\alpha,r) d\alpha. \]
As a last step -- just as before -- we only need to solve the equation $F'(r) = 0$. We get
\[ a = \frac{2}{\sqrt{2}+1} \cos \left(\frac{\pi}{16} \right)\,r. \]

\section{Construction}

In this section we show how to construct an approximation.
Let the radius $r$ of the approximated circle be given. The approximation will be given by the number
of used 1-steps and 2-steps. Depending on the approximation method one can compute
the sidelength ($a$) of the octagon that is given in the previous section (and it also can be
found in Table 4.)
By Table 2 we can compute the number of 1's and 2's used in the neighbourhood sequence
to obtain the desired octagon. The side length of the side $b$ can be computed by Table 1.

On Fig. 6 and 7 the number of 1's and 2's of the best approximating neighbourhood sequences are given
depending on the radius of the approximated circle.

\begin{figure}[ht]
  \centerline
  {\epsfig{file=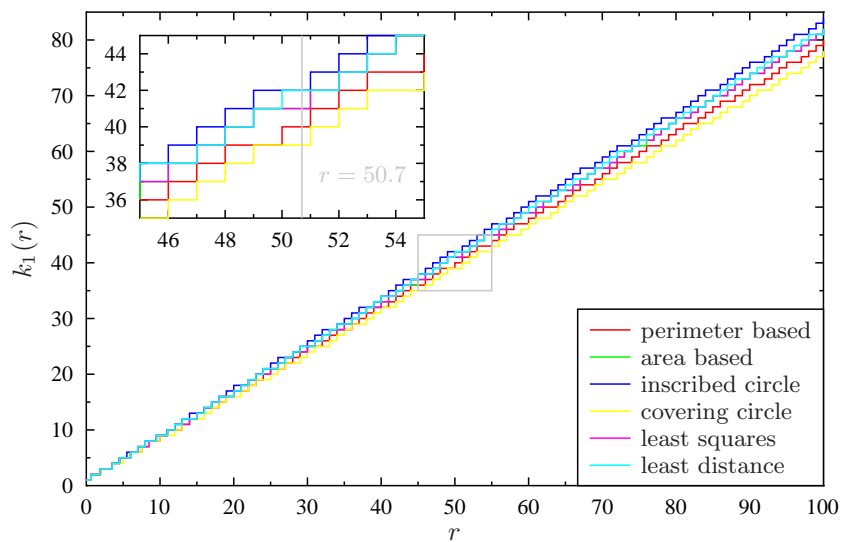} }
\caption{The number of 1's obtaining the best approximation}
\end{figure}

\begin{figure}[ht]
  \centerline
  {\epsfig{file=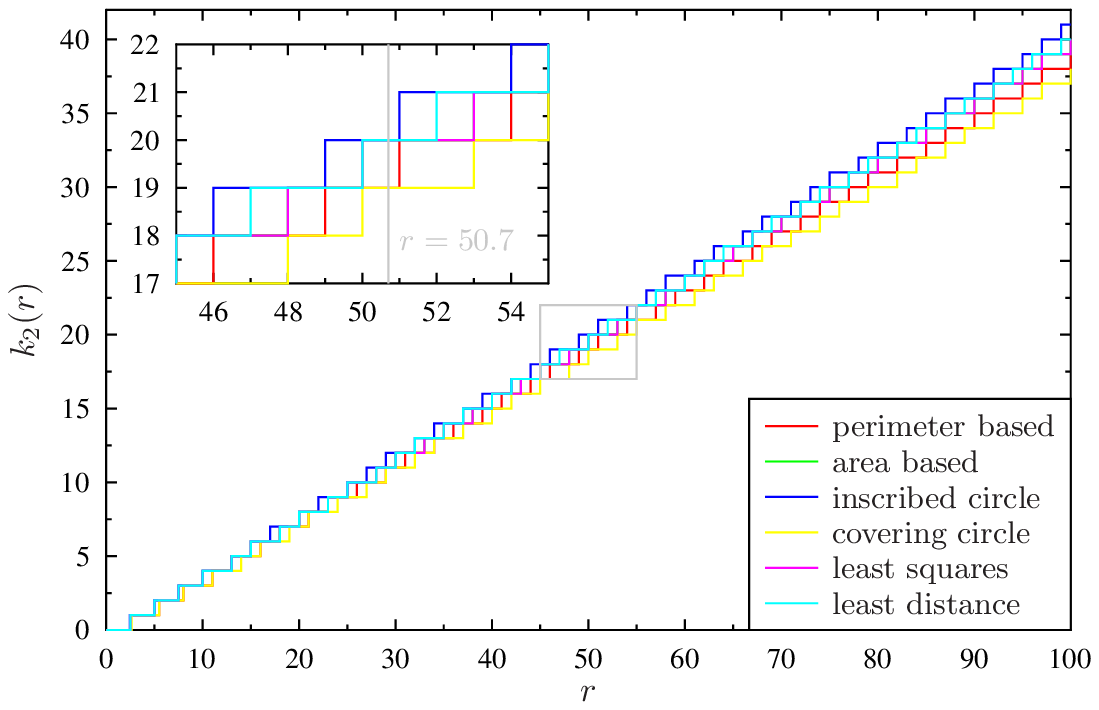} }
\caption{The number of 2's obtaining the best approximation}
\end{figure}

Now we show an example. Let $r=50.7$.
From the value of $a$ the values of $k_1$ and $k_2$ are computed, then from $k_1$ the value of $b$ can be
calculated.
  Table 3 shows the best approximations using the inner octagon descriptor.

\begin{table} \begin{center}
\begin{tabular}{|l|c|c|c|c|}
  \hline
                   & $a$ & $k_1$ & $k_2$ & $b$  \\                     \hline
   Perimeter based & 40 & 28 & 20 & 39.5980 \\
  Area based & 40 & 29 & 20 & 41.0122 \\
  Inscribed circle & 42 & 30 & 21 & 42.4264 \\
  Covering circle & 38 & 27 & 19 & 38.1838 \\
  Least squares & 40 & 29 & 20 & 41.0122 \\
  Least distance & 42 & 29 & 21 & 41.0122 \\
            \hline
\end{tabular} \end{center}
\begin{center}\emph{Table 3. The sidelengths $a$ and $b$ and the values $k_1, k_2$ approximating circle with
 $r=50.7$}
       \end{center}
\end{table}

In Fig. 8 the  best approximations obtained are shown.

\begin{figure}[ht]
  \centerline
  {\epsfig{file=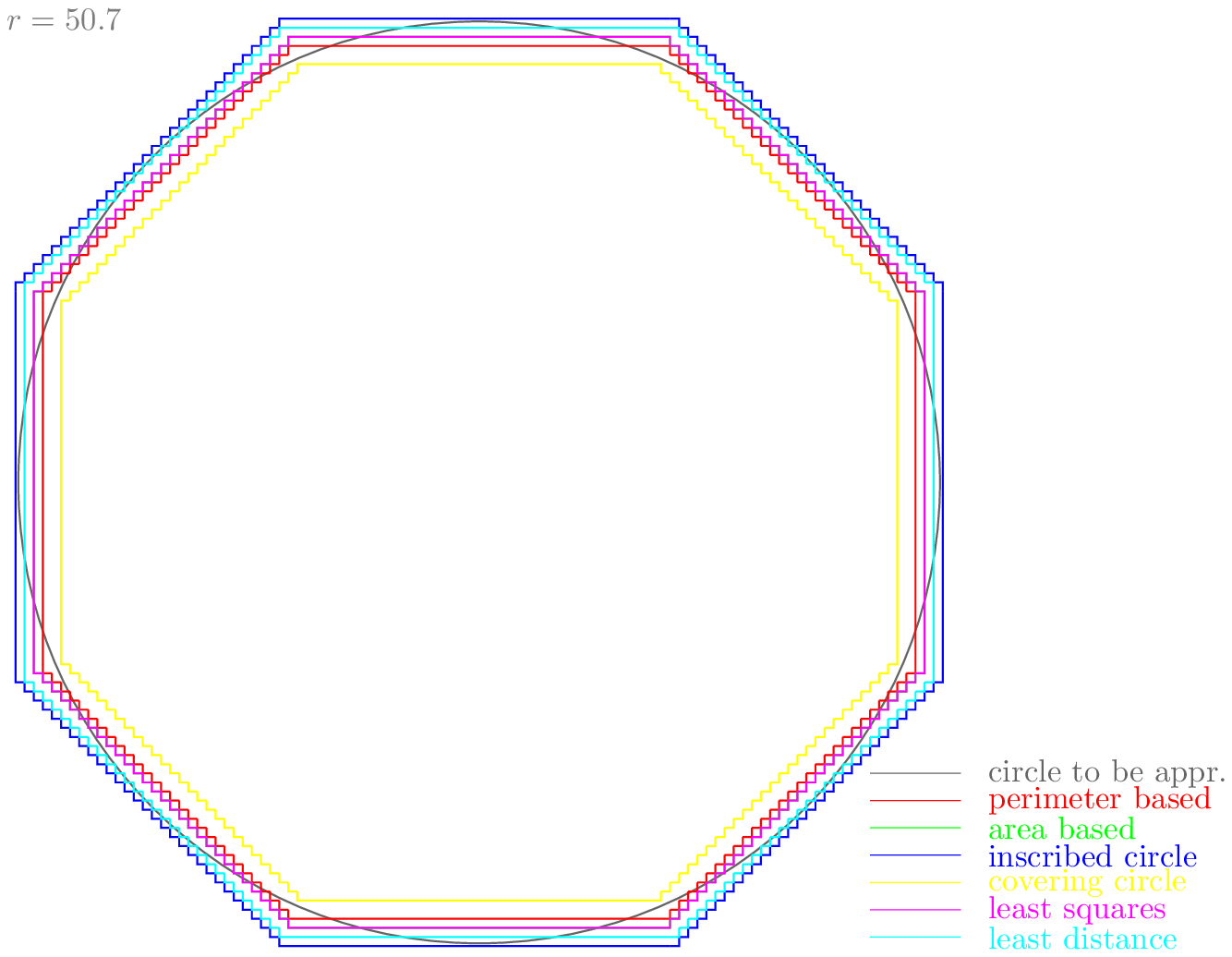,width=.94\hsize}}
\caption{The best approximating digital circles of the Euclidean circle $r=50.7$}
\end{figure}

\section{Comparison}
\par Here we summarize our main results in the form of a table. Table 4 shows the exact (continuous) sidelengths $a$ calculated by using the methods described in the preceding section. $r$ means the radius of the given circle.
With the help of Table 2 we can get the number of values 1 and 2 in the neighbourhood sequence we look for,
final formulae shown in Table 5.
Then $k = k_1 + k_2$ gives the minimum length of the sequence. After the $k$th element the sequence may contain arbitrary values. The order of the first $k$ elements is also arbitrary, since their permutation yields the same octagon. We also provide numerical results (rounded to 6 digits)
 to make the comparison of the coefficients easier. The developers only need to build Table 5 into their software supporting neighbourhood sequence based distance functions. In the table $[x]$ means the integer part of a real number $x$.
\bigskip
\par
\begin{table}
\begin{centering}
\begin{tabular}{c c c}
        \vspace{0.2cm}
    & \textbf{Exact sidelength}${\;(a)}$ & \textbf{Numerical sidelength} \\ \vspace{0.2cm}
    \textbf{Perimeter based}  & $\frac{\pi}{4}\,r$ & $0.785398\,r$\\ \vspace{0.2cm}
    \textbf{Area based} & $\sqrt{\frac{\pi}{2\left(1+\sqrt{2}\right)}}\,r$ & $0.806626\,r$\\ \vspace{0.2cm}
    \textbf{Inscribed circle} & $2\left(\sqrt{2}-1\right)\,r$ & $0.828427\,r$\\ \vspace{0.2cm}
    \textbf{Covering circle} & $\sqrt{2-\sqrt{2}}\,r$ & $0.765367\,r$ \\ \vspace{0.2cm}
    \textbf{Least squares} & $\frac{\pi}{4\left(\sqrt{2}+1\right)\ln \tan \left( \frac{5\pi}{16}\right)}\,r$ & $0.806852\,r$ \\ \vspace{0.2cm}
    \textbf{Least distance} & $\frac{2}{\sqrt{2}+1} \cos \left(\frac{\pi}{16} \right)\,r$ & $0.812509\,r$
\end{tabular}\\
\end{centering}
\begin{center} \emph{Table 4. The sidelength $a$ 
expressed by the radius $r$
in various approximations} \end{center}
\end{table}
\bigskip

\begin{table} \begin{center}
\begin{tabular}{|p{80pt}|p{80pt}p{80pt}|}
     \hline							 & \multicolumn{2}{c|}{Pixel based $(k_1\ge1)$} \\
     								 		 & \qquad \quad  $k_1$ & \qquad \quad  $k_2$ \\                     \hline
  Perimeter based  & $[0.785398r+1]$ & $[0.392699r-0.5]$  \\
  Area based       & $[0.806626r+1]$ & $[0.403313r-0.5]$  \\
  Inscribed circle & $[0.828427r+1]$ & $[0.414214r-0.5]$  \\
  Covering circle  & $[0.765367r+1]$ & $[0.382684r-0.5]$  \\
  Least squares    & $[0.806852r+1]$ & $[0.403426r-0.5]$  \\
  Least distance   & $[0.812509r+1]$ & $[0.406255r-0.5]$  \\
            \hline
\end{tabular}

\medskip

\begin{tabular}{|p{80pt}|p{80pt}p{80pt}|}
     \hline							 & \multicolumn{2}{c|}{Inner octagon} \\
     			               &   \qquad \ $k_1$ &   \qquad \ $k_2$ \\                     \hline
  Perimeter based  & $[0.555360r]$ & $[0.392699r]$ \\
  Area based       & $[0.570371r]$ & $[0.403313r]$ \\
  Inscribed circle & $[0.585786r]$ & $[0.414214r]$ \\
  Covering circle  & $[0.541196r]$ & $[0.382684r]$ \\
  Least squares    & $[0.570531r]$ & $[0.403426r]$ \\
  Least distance   & $[0.574531r]$ & $[0.406255r]$ \\
            \hline
\end{tabular}

\medskip

\begin{tabular}{|p{80pt}|p{80pt}p{80pt}|}
     \hline							 & \multicolumn{2}{c|}{Outer octagon} \\
     			               &   \qquad \ $k_1$ &  \qquad \quad \ $k_2$  \\                     \hline
  Perimeter based  & $[0.555360r]$ & $[0.392699r-0.5]$ \\
  Area based       & $[0.570371r]$ & $[0.403313r-0.5]$ \\
  Inscribed circle & $[0.585786r]$ & $[0.414214r-0.5]$ \\
  Covering circle  & $[0.541196r]$ & $[0.382684r-0.5]$ \\
  Least squares    & $[0.570531r]$ & $[0.403426r-0.5]$ \\
  Least distance   & $[0.574531r]$ & $[0.406255r-0.5]$ \\
            \hline
\end{tabular}
\end{center}
\begin{center}\emph{Table 5. The numerical values of $k_1, k_2$ 
 expressed by the radius $r$
in various approximations}
       \end{center}
\end{table}

\bigskip

\par We can choose the used descriptor and approximation method depending on our aim. For example, if we would like to have an exact number of pixels from a given distance of a center point, we shall use pixel descriptors with perimeter based approximation. If we would like to clusterize a plane with neighbourhood sequence generated discs, and we need clusters with fixed area, we shall use a convex hull descriptor (it makes no real difference whether we use the inner or the outer octagon) with area based approximation. The other approximations can be useful in visual applications, and their results can be more easily forecast, since we often work with the radius and the curve of the shapes when we imagine geometric ideas.

\section{Conclusions}
 The presented approximations work not only with integer radii
 (opposite to the previous papers, in which the best approximating neighbourhood sequence is computed
 to provide best approximation for the sequence of circles with integer radii). An example is shown for radius 50.7.
 Similarly, for any non-negative real radius
 one can compute the best approximation with our formulae.
   The digital discs and our formulae are useful for distance transforms, segmentation by colour clusterization and other algorithms in image processing. 
 A discussion of a three dimensional approximation is provided in \cite{6}, but the extension of most
 approaches of the problem to three and more dimensions is a matter of future research.
 There are some related results on the triangular grid (\cite{2,pre}) and in three dimensions on the
 face-centered and on the body-centered cubic grids (\cite{ro}).


\end{document}